\documentclass[12pt]{amsart}
\usepackage{amscd}
\def\NZQ{\Bbb}               
\def\NN{{\NZQ N}}

\def\frk{\frak}               

\def\pp{{\frk p}}

\def\Phi{{\frk n}}
\def\Phi{{\frk N}}
\def\opn#1#2{\def#1{\operatorname{#2}}} 
\opn\chara{char} \opn\length{\ell} \opn\pd{pd} \opn\rk{rk} \opn\projdim{proj\,dim} \opn\injdim{inj\,dim} \opn\rank{rank} \opn\depth{depth} \opn\grade{grade}
\opn\height{height} \opn\embdim{emb\,dim} \opn\codim{codim}

\opn\Tr{Tr} \opn\bigrank{big\,rank} \opn\superheight{superheight}\opn\lcm{lcm}
\opn\trdeg{tr\,deg}
\opn\reg{reg} \opn\lreg{lreg} \opn\ini{in} \opn\lpd{lpd}
\opn\size{size} \opn\Pf{Pf} \opn\GL{GL} \opn\SL{SL} \opn\mod{mod}
\opn\ord{ord} \opn\Gin{Gin}
\opn\Hilb{Hilb}\opn\adeg{adeg}\opn\std{std}\opn\ip{infpt}
\opn\div{div} \opn\Div{Div} \opn\cl{cl} \opn\Cl{Cl}
\opn\Spec{Spec} \opn\Supp{Supp} \opn\supp{supp} \opn\Sing{Sing} \opn\Ass{Ass} \opn\Min{Min}
\opn\Ann{Ann} \opn\Rad{Rad} \opn\Soc{Soc}
\opn\Syz{Syz} \opn\Im{Im} \opn\Ker{Ker} \opn\Coker{Coker} \opn\Am{Am} \opn\Hom{Hom} \opn\Tor{Tor} \opn\Ext{Ext} \opn\End{End} \opn\Aut{Aut} \opn\id{id}

\opn\nat{nat}
\opn\pff{pf}
\opn\Pf{Pf} \opn\GL{GL} \opn\SL{SL} \opn\mod{mod} \opn\ord{ord} \opn\Gin{Gin} \opn\Hilb{Hilb}
\opn\aff{aff} \opn\con{conv} \opn\relint{relint} \opn\st{st} \opn\lk{lk} \opn\cn{cn} \opn\core{core} \opn\vol{vol} \opn\link{link} \opn\star{star}
\opn\gr{gr}

\def\pot#1#2{#1[\kern-0.28ex[#2]\kern-0.28ex]}

\opn\dirlim{\underrightarrow{\lim}} \opn\inivlim{\underleftarrow{\lim}}
\let\union=\cup
\let\sect=\cap

\let\iso=\cong

\let\Sect=\bigcap

\let\to=\rightarrow

\def\Implies{\ifmmode\Longrightarrow \else
        \unskip${}\Longrightarrow{}$\ignorespaces\fi}
\def\implies{\ifmmode\Rightarrow \else
        \unskip${}\Rightarrow{}$\ignorespaces\fi}
\def\iff{\ifmmode\Longleftrightarrow \else
        \unskip${}\Longleftrightarrow{}$\ignorespaces\fi}

\let\:=\colon
\newtheorem{Theorem}{Theorem}[section]
\newtheorem{Lemma}[Theorem]{Lemma}
\newtheorem{Corollary}[Theorem]{Corollary}
\newtheorem{Proposition}[Theorem]{Proposition}

\newtheorem{Example}[Theorem]{Example}
\newtheorem{Examples}[Theorem]{Examples}

\let\epsilon\varepsilon
\let\phi=\varphi
\let\kappa=\varkappa
\def\qed{\ifhmode\textqed\fi
      \ifmmode\ifinner\quad\qedsymbol\else\dispqed\fi\fi}
\def\textqed{\unskip\nobreak\penalty50
       \hskip2em\hbox{}\nobreak\hfil\qedsymbol
       \parfillskip=0pt \finalhyphendemerits=0}
\def\dispqed{\rlap{\qquad\qedsymbol}}

\opn\dis{dis}
\def\pnt{{\raise0.5mm\hbox{\large\bf.}}}

\opn\Lex{Lex}

\begin{document}

\title{Criterions for  shellable multicomplexes}

\author{ Dorin Popescu}
\thanks{The  author was mainly supported by Marie Curie Intra-European
Fellowships MEIF-CT-2003-501046 and partially supported by the Ceres
program 4-131/2004 of the Romanian Ministery of Education and
Research}

\address{Dorin Popescu, Institute of Mathematics "Simion Stoilow", University of Bucharest,
P.O.Box 1-764, Bucharest 014700, Romania}
\email{dorin.popescu@imar.ro}

\maketitle
\begin{abstract}
After \cite{HP} the shellability of multicomplexes $\Gamma$ is given
in terms of some special faces of $\Gamma$ called facets. Here we
give a criterion for the shellability in terms of maximal facets.
Multigraded pretty clean filtration is the algebraic counterpart of
a shellable multicomplex. We give also a criterion for the existence
of a multigraded pretty clean filtration.
\end{abstract}

\section*{Introduction}

Cleanness is the algebraic counterpart of shellability for
simplicial complexes after \cite{D}. A kind of multigraded
"sequentially" cleaness the so called pretty cleaness was introduced
in \cite{HP}. Multigraded pretty cleaness implies sequentially
Cohen-Macaulay which remind us a well known result of Stanley
\cite{St} saying that shellable simplicial complexes are
sequentially Cohen-Macaulay. Pretty cleaness is the algebraic
counterpart of shellability of the so called multicomplexes (see
\cite{HP}). The aim of this paper is to find easy criterions for
multigraded pretty cleaness (see Theorem \ref{pret}) or for the
shellability of multicomplexes (see Theorem \ref{multi4}). The
Proposition \ref{facet} is important in the proof of
\cite[Proposition (10.1)]{HP}, where it was a consequence of some
results concerning standard pairs given in \cite{STV}. Here we give
an independent proof. Many useful examples are included. For
instance in \ref{dim} it is given a shellable multicomplex which has
a shelling $a_1,\ldots, a_r$ which does not satisfy the condition
$\dim S_1\geq \ldots \geq \dim S_r$ from \cite[Corollary (10.7)]{HP}
or here \ref{!} though certainly there is another shelling for which
it holds. Example \ref{bb} shows that there are shellable
multicomplexes which are not maximal shellable.

We express our  thanks to J. Herzog especially for some discussions
around Theorem \ref{multi4}.

\section{Preliminaries on pretty clean modules and multicomplexes}

 Let $R$ be a Noetherian ring, and $M$ a finitely generated
$R$-module. Then it is well known   that there exists a so called
{\sl prime filtration}
$$\mathcal{F}\: 0=M_0\subset M_1\subset \cdots \subset M_{r-1}
\subset M_r=M$$ that is such that
 $M_i/M_{i-1}\iso R/P_i$ for some   $P_i\in \Supp(M)$.
 We denote $\Supp(\mathcal{F})=\{P_1,\ldots, P_r\}$ and $r$ is
 called the length of $\mathcal{F}$.
It follows  that
\[
\Ass(M)\subset \Supp(\mathcal{F})\subset \Supp(M).
\]
If $\Supp(\mathcal{F})\subset \Min(M)$ then $\mathcal{F}$ is called
{\em clean}. $M$ is called {\em clean}, if $M$ admits a clean
filtration. $R$ is {\em clean} if it is a clean module over itself.
In particular, if $M$ is clean then $\Ass(M)= \Min(M)$. Let $\Delta$
be a simplicial complex on the vertex $\{1,\ldots,n\}$, $K$ a field
and $K[\Delta]$ the Stanley-Reisner ring.

\begin{Theorem}[Dress 1991]
\label{dr} $\Delta$  is a (non-pure) shellable simplicial complex if
and only if $K[\Delta]$  is a clean ring.
\end{Theorem}

A pure shellable simplicial complex is Cohen-Macaulay. So if $I$ is
a reduced monomial ideal of $S=K[x_1,\ldots,x_n]$ such that $S/I$ is
clean equidimensional then $S/I$ is Cohen-Macaulay. This result is
extended in \cite{HP} in a more general frame as we explain bellow.

 A prime filtration ${\mathcal F}\: 0=M_0\subset M_1\subset \ldots
\subset M_{r-1}\subset M_r=M$ of $M$ with $M_i/M_{i-1}=R/P_i$ is
called {\em pretty clean}, if for all $i<j$ for which $P_i\subset
P_j$ it follows that $P_i=P_j$. This means, roughly speaking, that a
proper inclusion $P_i\subset P_j$ is only possible if $i>j$. $M$ is
called {\em pretty clean}, if it has a pretty clean filtration. A
ring is called pretty clean if it is a pretty clean module over
itself. If ${\mathcal F}$ is pretty clean then  $\Supp({\mathcal
F})=\Ass(M)$.

\begin{Examples}{\em   Let $S=K[x,y]$ be the polynomial ring over the field $K$,
$I\subset S$ the ideal $I=(x^2,xy)$ and $R=S/I$. Then $R$ is pretty
clean but not clean. Indeed, $0\subset (x)\subset R$ is a pretty
clean filtration of $R$ with $(x)=R/(x,y)$,
 so that $P_1=(x,y)$ and $P_2=(x)$. $R$ is not clean since $\Ass(R)\neq
\Min(R)$. Note $R$ has a different prime filtration, namely,
${\mathcal G}: 0\subset (y)\subset (x,y)\subset R$ with factors
$(y)=R/(x)$ and $(x,y)/(y)=R/(x,y)$. Hence this filtration is not
pretty clean, even though $\Supp({\mathcal G})=\Ass(M)$. }
\end{Examples}

\begin{Example}
\label{ex1} {\em Let $R$ be a  UFD ring and $t_1,\ldots, t_s$ be
some irreducible elements in $R$, even equally some of them. Let $P$
be a prime ideal of $R$. Then $M:=R/I$, $I:=t_1\cdots t_s P$ is a
pretty clean module, even clean if $P$ does not contain any of
$(t_i)$. Indeed, consider the following filtration on $M$

$$M=M_{s+1}\supset M_s\supset \ldots \supset M_1\supset M_0=(0), $$
\noindent where $M_r:=(t_r\cdots t_s)/I$ for $1\leq r\leq s$. We
have
$$M_r/M_{r-1}\cong (t_r\cdots t_s)/(t_{r-1}\cdots t_s)\cong
R/(t_{r-1}).$$ \noindent and  $M_1\cong R/P$. A special type of this
example is given by
 $R=K[x,y]$ and $I=(x^2,xy)$.}
 \end{Example}

\begin{Proposition}[\cite{HP}]
\label{allthesame} Let $M$ be a pretty clean module. Then all pretty
clean filtrations of $M$ have the same length, namely their common
length equals

\noindent $\sum_{\pp\in\Ass(M)}$length$_M(H_{\pp}^0(M_{\pp}))$, that
is this number is bounded by the arithmetic degree of $M$, which is
$\sum_{\pp\in\Ass(M)}$length$_M(H_{\pp}^0(M_{\pp}))deg(R/\pp).$
\end{Proposition}

 $M$ is called {\em sequentially Cohen-Macaulay} if it a has
filtration
\[
0=C_0\subset C_1\subset C_2\subset\ldots\subset C_s=M
\]
such that $C_i/C_{i-1}$ is Cohen-Macaulay, and
$$\dim (C_1/C_0)<\dim(C_2/C_1)<\ldots <\dim(C_r/C_{r-1}).$$

\begin{Theorem}[Herzog, Popescu \cite{HP}]
\label{sequentially} Let $R$ be a local CM ring admitting a
canonical module $\omega_R$, and let $M$ be a pretty clean
$R$-module such that $R/P$ is Cohen-Macaulay for all $P\in \Ass M$.
Then $M$ is sequentially Cohen-Macaulay. Moreover, if $\dim R/P=\dim
M$ for all $P\in\Ass(M)$, then $M$ is clean and Cohen-Macaulay.
\end{Theorem}

Schenzel introduced in \cite{Sc} the so called {\em dimension
filtration}
\[
{\mathcal F}\: 0\subset D_0(M)\subset D_1(M)\subset \cdots \subset
D_{d-1}(M)\subset D_d(M)=M
\]
of $M$, which is defined by the property that $D_i(M)$ is the
largest submodule of $M$ with $\dim D_i(M)\leq i$ for
$i=0,\ldots,d=\dim M$.

\begin{Theorem}[Schenzel \cite{Sc}]
\label{seq} $M$ is sequentially CM, if and only if the factors in
the dimension filtration of $M$ are either 0 or CM.
\end{Theorem}

\begin{Theorem}[Herzog, Popescu \cite{HP}]
\label{sequentially} Let $R$ be a local CM ring admitting a
canonical module $\omega_R$, and let $M$ be a finitely generated
$R$-module  such that $R/P$ is Cohen-Macaulay for all $P\in \Ass M$.
Then $M$ is pretty clean if and only if $D_i(M)/D_{i-1}(M)$ is clean
for all $1\leq i\leq \dim M$.
\end{Theorem}

So in the above assumptions we may say that pretty clean means
sequentially clean. Now we pass to the multigraded case.
\begin{Proposition}[Herzog, Popescu \cite{HP}]
\label{tot} Let $I\subset S=K[x_1,\ldots,x_n]$ be a monomial ideal
such that $\Ass S/I$ is totally ordered by inclusion. Then $S/I$ is
pretty clean.
\end{Proposition}

A monomial ideal $I$ is called of {\em Borel type} if
$I:(x_1,\ldots,x_j)^{\infty}=I:x_j^{\infty}$ for all $1\leq j\leq
n$.
\begin{Corollary}[Herzog, Popescu, Vladoiu \cite{HPV}]
\label{bor} If the ideal $I\subset S$ is of Borel type then $S/I$ is
pretty clean and in particular sequentially CM.
\end{Corollary}

\noindent

Let $\Delta$ be a non-pure shellable simplicial complex on the
vertex set $\{1,\ldots,n\}$ and  $F_1,\ldots, F_r$ its shelling on
the facets of $\Delta$. For $i\geq 2$ we denote by $a_i$ the number
of facets of  $\langle F_1,\ldots, F_{i-1}\rangle\sect \langle
F_i\rangle$, and set $a_1=0$. Let $P_i=(\{x_j\}_{j\not \in F_i})$ be
the prime ideal associated to the facet $F_i$.

\begin{Proposition}[\cite{HP}]
\label{shelling numbers} The filtration $(0)=M_0\subset
M_1\subset\cdots\cdots M_{r-1}\subset M_r=K[\Delta]$ with
\[
M_i=\Sect_{j=1}^{r-i}P_{j}\quad \text{and}\quad M_i/M_{i-1}\iso
S/P_{r-i+1}(-a_{r-i+1})
\]
is a clean filtration of $S/I_{\Delta}$.
\end{Proposition}

Simplicial complexes correspond to the reduced monomial ideals of
$S$. What about general monomial ideals of $S$?

Let $\NN_\infty=\NN\union \{\infty\}$. For a subset
$\Gamma\subset\NN_\infty$ denote by ${\mathcal M}(\Gamma)$ the set
of all maximal elements of $\Gamma$. Let $a\in\Gamma$. Then $$\ip
a=\{i\:a(i)=\infty\}$$ is called the {\em infinite part} of $a$. A
subset $\Gamma\subset \NN^n_\infty$ is called a {\em multicomplex}
if
\begin{enumerate}
\item[(1)] for all $a\in\Gamma$ and all $b\in\NN^n_\infty$ with $b\leq a$ it follows that $b\in\Gamma$;

\item[(2)] for each $a\in\Gamma$ there exists $m\in{\mathcal M}(\Gamma)$ with $a\leq m$.
\end{enumerate}
The elements of a multicomplex are called {\em faces} and the
elements of ${\mathcal M}(\Gamma)$ are called {\em maximal facets}.
A face $a\in \Gamma$ is called {\em facet} if for any $m\in
{\mathcal M}(\Gamma)$ with $a\leq m$ it holds $\ip a=\ip m$.

Consider for example the multicomplex $\Gamma\in\NN^2_\infty$ with
faces $$\{a\: a\leq (0,\infty )\; \text{or}\; a\leq (2,0)\}.$$ Then
${\mathcal M}(\Gamma)=\{(0, \infty), (2,0)\}$ and ${\mathcal
F}(\Gamma)=\{(0,\infty), (2,0), (1,0)\}$. Besides its facets,
$\Gamma$ admits the infinitely many faces $(0,i)$ with $i\in \NN$.
\begin{Lemma}[\cite{HP}]
\label{finite} Each multicomplex has a finite number of facets.
\end{Lemma}
Let $\Gamma$ be a multicomplex, and  let $I(\Gamma)$ be the
$K$-subspace in $S$ spanned by all monomials $x^a$ such that
$a\not\in \Gamma$. This is an ideal in $S$ and the correspondence
$\Gamma\to I(\Gamma)$ gives a bijection between  the multicomplexes
$\Gamma$ of $\NN^n_\infty$ and the monomial ideals $I$ of $S$.

 Let  $a\in \NN^n$, $m\in \NN^n_\infty$ with $m(i)\in\{0,\infty\}$
 and $\Gamma(m)$ the multicomplex generated by $m$, that is the set
 of all $u\in \NN^n_\infty$ with $u\leq m$. Actually given $m_1.\ldots,m_r\in \NN^n$ we denote by $\Gamma(m_1,\ldots,m_r)$
 the set of all $u\in \NN^n$ with $u\leq m_i$ for some $i$.  The sets of the form $S=a+S^*$,
 where $S^*=\Gamma(m)$, are called {\em Stanley
 sets}. The {\sl dimension of S} is defined to be the dimension of
 $S^*$.
 A multicomplex $\Gamma$ is {\em shellable} if the facets of
$\Gamma$ can be ordered $a_1,\ldots, a_r$ such that
\begin{enumerate}
\item[(1)] $S_i= \Gamma(a_i)\setminus\Gamma(a_1,\ldots, a_{i-1})$ is a Stanley set for $i=1,\ldots,r$, and

\item[(2)] whenever $S_i^*\subset S_j^*$, then $S_i^*=S_j^*$ or $i>j$.
\end{enumerate}

\begin{Theorem}[Herzog, Popescu \cite{HP}]
\label{multi2} The multicomplex $\Gamma$ is shellable if and only if
$S/I(\Gamma)$ is a multigraded pretty clean ring.
\end{Theorem}

 Note that in the above example the shelling could be
 $\{(0,\infty),(1,0),(2,0)\}$ and so the first Stanley set is given
 by the axe $\{(0,s):s\in \NN\}$ and the second Stanley set, respectively the third
 one are the points $(1,0)$, $(2,0)$.

\begin{Corollary}[\cite{HP}]
\label{!} A multicomplex $\Gamma$ is shellable if and only if there
exists an order $a_1,\ldots,a_r$ of the facests such that for
$i=1,\ldots,r$ the sets $S_i=\Gamma(a_i)\setminus
\Gamma(a_1,\ldots,a_{i-1})$ are Stanley sets with $\dim S_1\geq \dim
S_2\geq\ldots \geq\dim S_r$.
\end{Corollary}

\section{Pretty clean modules}

Let $R$ be a Noetherian ring, and $M$ a finitely generated
$R$-module.

Let
$${\mathcal M}: M_0=0\subset M_1\subset \ldots \subset M_s=M$$
\noindent be a filtration of $M$. Inspired by \cite{D} we call
${\mathcal M}$
 {\sl almost clean} if for all $1\leq i\leq s$ there exists
 a prime ideal $P\in \Ass_R(M) $ such that $P=(M_{i-1}:x)$ for all
 $x\in M_i\setminus M_{i-1}$.

 \begin{Lemma}
 \label{almost} Every finitely generated $R$-module $M$ has an
 almost  clean filtration.
 \end{Lemma}

 \begin{proof} Let $\Ass_R(M)=\{P_1,\ldots,P_t \}$ and $(0)=
 \cap_{i=1}^t N_{i}$ be an irredundant primary decomposition
 of $(0)$ in $M$, where $N_i$ is a $P_i$-primary submodule of $M$.
  We may suppose the notation such that for all $i<j$
such that $P_i\subset P_j$ it follows $P_i=P_j$.
   Set $U_j=\cap_{i=1}^j N_{i}$. We get a
 filtration

 $$(0)=U_t\subset \ldots \subset U_1\subset U_0=M$$

 \noindent such that $U_{i-1}/U_i\cong U_{i-1}+N_{i}/ N_{i} \subset
M/N_{i}$. But $U_{i-1}\not = U_i$ because the primary decomposition
is irredundant. Then
$\{P_i\}=\Ass_R(M/N_{i})=\Ass_R(U_{i-1}/U_i)\not=\emptyset$ and so
$\Ass_R(U_{i-1}/U_i)=\{P_i\}$. Thus we reduce to the case
$T=U_{i-1}/U_i$. In this case set $V_k=P_i^kT_{P_i}\cap T$ for
$k\geq 0$. We get a filtration
$$(0)\subset \ldots \subset V_1\subset V_0=T$$
\noindent such that there exists an injection $V_{j-1}/V_j\to
P_i^{j-1} T_{P_i}/P_i^j T_{P_i}$, the last module being a linear
space over the fraction field of $R/P_i$. Thus $V_{j-1}/V_j$ is
torsionless over $R/P_i$, which is enough.

\end{proof}

\begin{Corollary}
\label{reg} If $R/P$ is regular of dimension $\leq 1$ for all $P\in
\Ass_R(M)$ then $M$ is pretty clean.
\end{Corollary}

For the proof note only that the quotients $V_{j-1}/V_j$ from the
above proof are in this case free and so clean.

Next proposition is an extension of \cite[Proposition (5.1)]{HP}.
 Let $S=K[x_1,\ldots,x_n]$ and $I\subset S$ a monomial ideal.
For every $P\in \Ass(S/I)$ let $J_P\subset \{1,\ldots,n\}$ be the
subset of all $j$ such that $x_j\in P$. Clearly $P$ is monomial and
so $\height(P)=|J_P|$ and $\dim P=n-|J_P|$.

\begin{Theorem}
\label{pret} Suppose that for all integer $d>0$ such that
$\Ass^d(S/I)\not =\emptyset$ it holds $|\bigcup_{P\in \Ass^d(S/I)}
J_P|\leq n-d+1$. Then $S/I$ is pretty clean.

\end{Theorem}

\begin{proof}
We follow the proof of Proposition \ref{tot} given in \cite{HP}. Let
$I=\bigcap_{P\in \Ass(S/I)} Q_P$ be the irredundant primary monomial
decomposition of $I$, where $Q_P$ is $P$-primary. Set $d_P=\dim P$
and for $d_1<\ldots <d_r$ be the integers which appear really among
$(d_P)$. Set $U_i=\bigcap_{P\in \cup_{j>d_i}\Ass^j(S/I)}\ Q_P$. By
\cite{Sc} the dimension filtration is given by $D_{d_i}(S/I)=U_i/I$.
Using Theorem \ref{sequentially} it is enough to show that
$U_i/U_{i-1}$ is clean. Let $S'$ be the polynomial ring over $K$ in
the variables $x_j$ with $j\in\bigcup_{P\in \Ass^{d_i}(S/I)}$ and
set $P'=P\cap S'$, $Q_{P'}=Q_P\cap S'$ for $P\in \Ass^{d_i}(S/I)$.
We have $P=P'S$, $Q_P=Q_{P'}S$ and $U_i=U'_iS$ for $U'_i=U_i\cap
S'$.. But $U'_i/U'_{i-1}$ is clean by Corollary \ref{reg} since
$\dim S' =n-d_i+1$ by assumption (that is $\dim S/P\leq 1$ for all
$P\in \Ass^{d_i}(S/I)$). Then by base change $U_I/U_{i-1}\iso
U'_i/U'_{i-1}\otimes_{S'}S$ is a clean module.

\end{proof}

\begin{Example}
\label{?} {\em Let
$I=(x_1^2,x_1x_2^2x_3,x_1x_3^2,x_2^2x_4^2,x_2x_3^2x_4)\subset
S=K[x_1,\ldots,x_4]$. We have
$$I=(x_1,x_2)\cap (x_1,x_4)\cap
(x_1^2,x_2^2,x_3^2)\cap (x_1^2,x_1x_3,x_3^2,x_4^2)$$ and so
$\Ass(S/I)=\{(x_1,x_2),(x_1,x_4),(x_1,x_2,x_3),(x_1,x_3,x_4)\}$. An
algorithm to find a monomial primary decomposition is given in
\cite{Vi}. Then $S/I$ is pretty clean by Theorem \ref{pret}.}
\end{Example}

In \cite{HP} is given an example of $R$-module $M$  which is not
pretty clean but has a prime filtration $\mathcal F$ with
$\Supp({\mathcal F})=\Ass_R(M)$. The following shows that there are
$R$-modules $M$ for which there exist no prime filtration ${\mathcal
F}$ with $\Supp({\mathcal F})=\Ass_R(M)$. Modules which has a
filtration $\mathcal F$ with $\Supp({\mathcal F})=\Ass_R(M)$ are
studied in different papers (see e.g. \cite{Li}) after Eisenbud's
question from \cite{E}.

\begin{Example} {\em Let $S=K[x_1,\ldots,x_4]$, $I=P_1\cap P_2$ for
$P_1 =(x_1,x_2)$, $P_2=(x_3,x_4)$ and $M=S/I$. Suppose that there
exists a filtration $\mathcal F$ of $M$ such that $\Supp({\mathcal
F})=\Ass_S(M)$. Then $\mathcal F$ is a clean filtration of $M$
because $\Ass_S(M)=\Min(M)=\{P_1,P_2\}$. Note that $P_i$ could
appear in $\mathcal F$ only of $1=\length_{S_{P_i}}(M_{P_i})$-times.
Thus $\mathcal F$ has the form $(0)\subset N\subset M$, where we may
suppose that $M/N\iso S/P_1$ and $N\iso S/P_2$. But this is not
possible because then $N=P_1/I$ is not cyclic. Thus $M$ has no
filtration $\mathcal F$ with $\Supp(\mathcal F)=\Ass_S(M)$. Note
that the hypothesis of Theorem \ref{pret} do not hold in this
frame.}
\end{Example}

\section{Multicomplexes}

The following proposition is stated in the proof of
\cite[Proposition (10.1)]{HP} using the standard pairs of
\cite{STV}. We think that this result deserves a direct proof which
we give bellow.
\begin{Proposition}
\label{facet} Let $\Gamma $ be a multicomplex. The arithmetic degree
of $S/I(\Gamma)$ is exactly the number of the facets ${\mathcal
F}(\Gamma)$ of $\Gamma$.
\end{Proposition}
\begin{proof}
Let $\phi$ be the map from $\Gamma$ to the monomial $K$-basis of
$S/I(\Gamma)$, given by $u\rightarrow x^{\tilde u}$, where ${\tilde
u}(j)=u(j)$ if $j\not \in \ip(u)$ and ${\tilde u}(j)=0$ otherwise.
Let $P\in \Ass(S/I(\Gamma))$ and ${\mathcal F}(\Gamma,P)=\{u\in
{\mathcal F}(\Gamma):P_u=P\}$. We claim that the restriction of
$\phi$ to ${\mathcal F}(\Gamma,P)$ is injective. Indeed, suppose
that we have ${\tilde u}={\tilde v}$ for some $u,v\in {\mathcal
F}(\Gamma ,P)$. We have $\ip(u)=\ip(v)$ since $P_u=P_v=P$ and so it
follows $u=v$.

Now let $u<v$ be two faces of $\Gamma$. Then $\ip(u)\subset \ip(v)$
and we have equality if and only if $(I(\Gamma(v)):x^{\tilde u})$ is
a $P_u$-primary ideal. Indeed if $\ip(u)=\ip(v)$ then
$(I(\Gamma(v)):x^{\tilde u})=(\{x_k^{{\tilde v}(k)-{\tilde u}(k)+1}:
k\not\in \ip(u)\})$, that is a $P_u$-primary ideal. Conversely, for
every $k\not\in \ip(u)$ suppose that a power $x_k^{\alpha_k}\in
(I(\Gamma(v)):x^{\tilde u})$. Thus $x_k^{\alpha_k} x^{\tilde u}\in
I(\Gamma(v))$, that is ${\tilde u}+\alpha_k \epsilon_k\not\in
\Gamma(v)$, where $\epsilon_k$ is the $k$-unitary vector. Then
$u(k)+\alpha_k>v(k)$, that is $k\not\in\ip(v)$.

We claim that given a face $u\in\Gamma$ it holds $u\in {\mathcal
F}(\Gamma)$ if and only if $\phi(u)\in
H_{P_uS_{P_u}}^0((S/I(\Gamma))_{P_u})$. Indeed, let $u$ be a facet,
then for each maximal facet $v$ with $u\leq v$ it holds
$\ip(u)=\ip(v)$. This happens if $(I(\Gamma(v)):x^{\tilde u})$ is a
$P_u$-primary ideal for any maximal facet $v\geq u$. It follows
$$\Sect_{v\in {\mathcal M}(\Gamma), u\leq v}(I(\Gamma(v)):x^{\tilde
u})$$ is a $P_u$-primary ideal. The converse is also true because
$$P_u=\sqrt{(\Sect_{v\in {\mathcal M}(\Gamma), u\leq v}(I(\Gamma(v)):x^{\tilde
u})}\subset \sqrt{I(\Gamma(v)):x^{\tilde u})}=P_v,$$ that is
$\ip(u)=\ip(v)$. So $u$ is a facet if and only if
$(I(\Gamma)S_{P_u}:x^{\tilde u})$ is $P_uS_{P_u}$-primary ideal and
it follows that there exist exactly
$dim(H_{PS_{P}}^0((S/I(\Gamma))_{P})$-facets $u$ in $\Gamma$ with
$P_u=P$.

\end{proof}

Let $\Gamma\subset \NN_{\infty}^n$ be a multicomplex and $a,b\in
\Gamma$. We call $a$ a {\it lower neighbour} of $b$ if there exists
an integer $k$, $1\leq k\leq n$ such that

\begin{enumerate}

\item[(i)] $a(i)=b(i)$ for $i\neq k$,

\item[(ii)] either  $a(k)+1=b(k)<\infty$, or $a(k)<\infty$ and $b(k)=\infty$.
\end{enumerate}

\medskip
It is easy to see that the multicomplex $\Gamma$ is  {\it shellable}
if the facets of $\Gamma$ can be ordered $a_1,\ldots, a_r$ such that
\begin{enumerate}
\item{} $a_1\in \{0,\infty\}^n,$

 \item{} for $i=2,\ldots,r$ the maximal
facets of $\langle a_1,\ldots,a_{i-1}\rangle\cap\langle a_i \rangle$  are lower neighbours of $a_i$;

\item{} for each $k\not\in \supp\ a_i$, $1\leq k\leq n$ such that $a_i(k)>0$ there exists a maximal facet $w$ of $\Gamma(a_1,\ldots,a_{i-1})\cap\Gamma(a_i)$
such that $w(k)<a_i(k).$

\item{} for all $1\leq j<i\leq r$ such that  $\supp a_j\subset
 \supp a_i$, it follows that $\supp a_j=\supp a_i$.
\end{enumerate}
Actually $\Gamma$ satisfies (a) above for any $i>1$  if and only if
it satisfies (2),(3) above, and  it satisfies (a) for $i=1$ if and
only if (1) holds. Also $\Gamma$ satisfies (b) above if and only if
it satisfies (4) above. There are orders of the facets of some
shellable multicomplexes which satisfy (1)-(3) but not (4) as shows
the following:

\begin{Example}
\label{(4)}{\em Let $a=(\infty,0,\infty,\infty)$,
$b=(1,1,\infty,0)$,$c=(0,2,\infty,\infty)$ and $\Gamma=\langle
a,b,c\rangle$. Then ${\mathcal
F}(\Gamma)=\{a,b,c,(o,1,\infty,\infty)\}$. We may order these facets
in the following way $u_1=a$,$u_2=(0,1,\infty,\infty)$, $u_3=b$,
$u_4=c$. Note that $\langle u_1\rangle \cap\langle u_2\rangle$ has
just one maximal facet $(0,0,\infty,\infty)$ which is a neighbour of
$u_2$. Also $\langle u_1,u_2\rangle \cap \langle u_3\rangle$ has two
maximal facets $(1,0,\infty,0)$, $(0,1,\infty,0)$, both being
neighbours of $u_3$. Finally note that $\langle
u_1,u_2,u_3\rangle\cap \langle u_4\rangle$ has just one maximal
facet $u_2$ which is a neighbour of $u_4$. So it is easy to see that
this order satisfies (1)-(3), but not (4) because
$P_{u_3}=(x_1,x_2,x_4)\supset (x_1,x_2)=P_{u_4}$. Actually, $\Gamma$
is shellable because of Theorem \ref{pret} or \cite[Proposition
5.1]{HP}.}
\end{Example}

Also there are shellable multicomplexes $\Gamma$ which have
shellings $u_1,\ldots u_r$ for which the Stanley sets
 $S_i=\Gamma(a_i)\setminus\Gamma(a_1,\ldots, a_{i-1})$ does not satisfy
  $\dim S_1\geq \dim S_2\geq \ldots\geq \dim S_r$ as shows the
  following:

  \begin{Example}
  \label{dim}{\em Let $a=(0,\infty,1,\infty)$, $b=(0,0,2,\infty)$,
  $c=(\infty,\infty,1,0)$ and $\Gamma=\langle a,b,c\rangle$.
  $\Gamma$ has apart of the maximal faces $a,b,c$ and the following
  facets: $d=(\infty,\infty,0,0)$, $e=(0,\infty,0,\infty)$. Choose
  the order $u_1=d$, $u_2=e$, $u_3=a$, $u_4=b$, $u_5=c$. Note that
$\langle u_1\rangle \cap\langle u_2\rangle$ has just one maximal
facet $(0,\infty,0,0)$ which is a neighbour of $u_2$ and $\langle
u_1,u_2\rangle \cap \langle u_3\rangle$ has one maximal facet
$(0,\infty,0,\infty)$ a neighbour of $u_3$. Also note that $\langle
u_1,u_2,u_3\rangle\cap \langle u_4\rangle$ has just one maximal
facet $(0,0,1,\infty)$ which is a neighbour of $u_4$ and $\langle
u_1,u_2,u_3,u_4\rangle\cap \langle u_5\rangle$ has two maximal
facets $(\infty,\infty,0,0)$, $(0,\infty,1,0)$ both being neighbours
of $u_5$. We have $P_{u_1}=(x_3,x_4)=P_{u_5}$,
$P_{u_2}=(x_1,x_3)=P_{u_3}$, $P_{u_4}=(x_1,x_2,x_3)$. Since $\dim
P_{u_5}>\dim P_{u_4}$ the filtration does not satisfy the dimension
condition above though it is pretty clean.}
\end{Example}

Let $\Gamma\subset \NN_{\infty}^n$ be a multicomplex and $u_1,\ldots,u_r$
 its maximal facets so $\Gamma=\langle u_1,\ldots,u_r\rangle$.

\begin{Lemma}
\label{basechange} If $\ip u_j=\ip u_1$ for all $j\geq 1$, that is
$I(\Gamma)$ is primary ideal, then $\Gamma$ is shellable and
$S/I(\Gamma)$ is clean.
\end{Lemma}

\begin{proof} This lemma is a consequence of Theorem
\ref{pret} or of \cite[Proposition 5.1]{HP} but we prefer to give
here the proof since it is elementary.  Note that a monomial primary
ideal $Q$ can be seen as the extension of a primary ideal $Q'$
associated to a maximal ideal in a polynomial ring $S'$ in fewer
variables which enter really in the generators of $Q$. Then $S'/Q'$
is a clean module and so by base change $S/Q$ is too.
\end{proof}

$\Gamma$ is {\it maximal shellable} if the maximal facets of $\Gamma$
can be ordered $u_1,\ldots, u_r$ and there exists $s$, $1\leq s\leq r$ such that
\begin{enumerate}
\item{}  $\ip u_1=\ip u_j,$ for all $1\leq j\leq s$

\item{} for $i=s+1,\ldots,r$ the maximal facets of $\langle u_1,\ldots,u_{i-1}\rangle\cap \langle u_i\rangle$
 differ from $u_i,$
only in one component,

\item{} for all $s\leq j<i\leq r$ such that $\ip\ u_j\subset \ip\ u_i$ it follows $\ip\ u_j=\ip\ u_i$.
\end{enumerate}

Suppose $\Gamma$ satisfies the above conditions, so it is maximal
shellable. Fix an $i=s+1,\ldots r$ and let $w_{i1},\ldots,w_{ic}$ be
the maximal facets of $\Gamma(u_1,\ldots,u_{i-1})\cap\Gamma(u_i)$.
Thus for each $1\leq j\leq c$ there exists just one $\lambda_j$,
$1\leq\lambda_j\leq n$ such that $w_{ij}(\lambda_j)<u_i(\lambda_j)$
and so $w_{ij}(\lambda_j)\in \NN$.  Set
$f_i=\Pi_{j=1,w_{ij}(\lambda_j)<\infty}^c\
 x^{w_{ij}(\lambda_j)+1}_{\lambda_j}.$ We claim that

\[
\Sect_{s=1}^{i-1} I(\Gamma(u_s))+I(\Gamma(u_i)) =I(\Gamma(u_i))+(f_i).
\]
The monomial $f_i$ has the form $f_i=x^{t_i}$ for some $t_i\in
\NN^n$. By definition of $f_i$ we have $t_i\leq u_i$, that is
$t_i\in \Gamma(u_i)$ and $t_i\not\in \Gamma(w_{ij})$ for all $1\leq
j\leq c$. Thus $t_i\in
\Gamma(u_i)\setminus\Gamma(w_{i1},\ldots,w_{ic})$. Since
$\Gamma(w_{i1},\ldots,w_{ic})=\Gamma(u_1,\ldots,u_{i-1})\cap
\Gamma(u_i)$ it follows $t_i\not\in \Gamma(u_1,\ldots,u_{i-1})$,
that is

$$f_i=x^{t_i}\in I(\Gamma(u_1,\ldots,u_{i-1}))=\Sect_{s=1}^{i-1}
I(\Gamma(u_s)).$$

\noindent Conversely,  let $x^q\in
I(\Gamma(u_1,\ldots,u_{i-1}))\setminus I(\Gamma(u_i))$, that is
$q\not\in \Gamma(u_1,\ldots,u_{i-1})$ and $q\in \Gamma(u_i)$. Then
$q\not \in \Gamma(w_{i1},\ldots,w_{ic})=\cup_{j=1}^c \Gamma(w_{ij})$
and $q\leq u_i$. Thus $u_i(k)\geq q(k)>w_{ij}(k)$ for at least one
$k$ and so $k=\lambda_j$ and it follows that $x^q\in (f_i)$.

Set $a_i=deg\ f_i$ for $i>1$ and $a_1=0$. We obtain the following isomorphisms of grade $S$-modules

\[
(\Sect_{s=1}^{i-1} I(\Gamma(u_s)))/(\Sect_{s=1}^{i} I(\Gamma(u_s)))\iso (\Sect_{s=1}^{i-1} I(\Gamma(u_s)))+ I(\Gamma(u_i))) /I(\Gamma(u_i)) =
\]
\[
(I(\Gamma(u_i))+(f_i))/I(\Gamma(u_i))\iso (f_i)/f_iI(\Gamma(u_i))\iso S/(I(\Gamma(u_i)):f_i)(-a_i).
\]

By construction of $f_i$ we see that
$x_{\lambda_j}^{u_i(\lambda_j)-w_{ij}(\lambda_j)}f_i\in
I(\Gamma(u_i))$ if $\lambda_j\not \in \ip u_i$. If $0\leq
u_i(k)<\infty$ and $k$ is not a $\lambda_j$ then $x_k$ does not
enter in $f_i$ and so $x_k$ enters in $(I(\Gamma(u_i)):f_i)$ only at
the power he had in $I(\Gamma(u_i))$. However we see that
$(I(\Gamma):f_i)$ is a irreducible $P_{u_i}$-primary ideal.
 Finally note that the condition (3) says that for all $i>j$ such
that $P_{u_i}\subset P_{u_j}$ it follows $P_{u_i}=P_{u_j}$. Thus we have shown:

\begin{Proposition}
\label{multi3} Let $\Gamma\subset \NN_{\infty}^n$ be a maximal shellable multicomplex and $u_1,\ldots,u_r$ its shelling. Then there exists a filtration of
$S/I(\Gamma)$

$(0)=M_0\subset M_1\subset\cdots\cdots M_{r-s}\subset M_{r-s+1}=S/I(\Gamma)$ with
\[
M_i=\Sect_{j=1}^{r-i}I(\Gamma(u_j))\quad \text{and}\quad M_i/M_{i-1}\iso S/J_i(-a_{r-i+1})
\]
for some irreducible ideals $J_i$ associated to $P_{u_{r+1-i}}$,
 for $i\leq r-s$ and the primary ideal $J_{r-s+1}=I(\Gamma(u_1,\ldots,u_s))$.
\end{Proposition}

\begin{Theorem}
\label{multi4} If the  multicomplex $\Gamma\subset \NN_{\infty}^n$
is maximal shellable then $\Gamma$ is shellable in particular
$S/I(\Gamma)$ is a pretty clean ring.

\end{Theorem}

\begin{proof}
By the above proposition it is enough to see that that the factors from the above filtration are clean. This follows by the Lemma \ref{basechange}.

\end{proof}

We end this section with some examples.

\begin{Example}{\em Let $J=(x_1^2,x_2^2,x_3,x_4)\cap
(x_1,x_2,x_3^2,x_4^2)$ and $T= (x_1,x_2,x_5^2,x_6^2)$
 be primary
ideals  in $S=K[x_1,\ldots,x_6]$. Set $I=J\cap T$. Then $S/I$ is not Cohen-Macaulay since from the exact sequence

$$0\to S/I\to S/J\oplus S/T\to S/J+T\to 0$$

\noindent we get $\depth (S/I)=1$. It follows that $S/I$ is not pretty clean too and so $\Gamma(I)$ cannot be maximal shellable. This indeed is the case since
one can take $s=2,r=3$, $u_1=(1,1,0,0,\infty,\infty)$, $u_2=(0,0,1,1,\infty,\infty)$ and $u_3=(0,0,\infty,\infty,1,1)$. For this order one can see that
$\Gamma(u_1,u_2,u_3)$ is not maximal shellable since the condition (2) does not hold.}
\end{Example}

The following example shows that there exist shellable
multicomplexes which are not maximal shellable.

\begin{Example}
\label{bb}{\em Let $J=(x_1^2,x_2^2,x_3,x_4)\cap
(x_1,x_2,x_3^2,x_4^2)$ and $L= (x_1^2,x_2,x_3,x_5^2)$
 be primary
ideals  in $S=K[x_1,\ldots,x_5]$. Set $I=J\cap L$. Then $S/I$ is
Cohen-Macaulay but  $\Gamma(I)$ is not maximal shellable. This
indeed is the case since one can take $s=2,r=3$,
$u_1=(1,1,0,0,\infty)$, $u_2=(0,0,1,1,\infty)$ and
$u_3=(1,0,0,\infty,1)$. The maximal facets of
$\Gamma(u_1,u_2)\cap\Gamma(u_3)$ are $w_1=(1,0,0,0,1)$ and $w_2=
(0,0,0,1,1)$. Clearly $w_1$   satisfies the condition (2) but $w_2$
not. Also note that $J=I+Sx_4^2$ and so the factors of the
filtration $(0)\subset J/I\subset S/I$ are all cyclic given by
primary ideals. By Lemma \ref{basechange} we see that $S/I$ is
pretty clean and so shellable. }
\end{Example}

\begin{Example}{\em Let $J=(x_1^2,x_2^2,x_3,x_4)\cap
(x_1,x_2,x_3^2,x_4^2)$ and $T= (x_1,x_2,x_3,x_5^2)$
 be primary
ideals  in $S=K[x_1,\ldots,x_5]$. Set $I=J\cap T$. Then $S/I$ is Cohen-Macaulay because  $\Gamma(I)$ is  maximal shellable. This indeed is the case since one
can take $s=2,r=3$, $u_1=(1,1,0,0,\infty)$, $u_2=(0,0,1,1,\infty)$ and $u_3=(0,0,0,\infty,1)$. The only maximal facet of $\Gamma(u_1,u_2)\cap\Gamma(u_3)$ is
$w=(0,0,0,1,1)$. Clearly $w$ satisfies the condition (2).}
\end{Example}


\begin{thebibliography}{1}


\bibitem{BH}  W.\ Bruns,\ J. Herzog, {\sl Cohen-Macaulay rings}, Revised
Edition,  Cambridge, 1996.

\bibitem{D} A.\ Dress, A new algebraic criterion for shellability, Beitrage
zur Alg. und Geom., {\bf 340}(1), (1993),45--55.

\bibitem{E} D.\ Eisenbud, {\sl Commutative algebra , with a view toward
geometry}, Graduate Texts Math. Springer, 1995, page 93.

\bibitem{HP} J.\ Herzog, D.\ Popescu, Finite filtrations of modules
and shellable multicomplexes,  Preprint IMAR no 4/2005, Bucharest,
2005.

\bibitem{HPV}  J.\ Herzog, D.\ Popescu, M.\ Vladoiu, On the Ext-modules of
ideals of Borel type, Contemporary Math. {\bf 331} (2003), 171-186.

\bibitem{Li} A.\ Li, Associated prime filtrations of finitely
generated modules over noetherian rings, Communications in Alg.,
{\bf 23(4)}, (1995), 1511-1526.


\bibitem{Sc} P.\ Schenzel, On the dimension filtration and
Cohen-Macaulay filtered modules, Proceed. of the Ferrara meeting in
honour of Mario Fiorentini, ed. F. Van Oystaeyen, Marcel Dekker,
New-York, 1999.



\bibitem{St} R.\ P.\ Stanley, {\sl Combinatorics and Commutative Algebra},
Birkh\"auser, 1983.



\bibitem{STV} B.\ Sturmfels, N.\ V.\ Trung, W. Vogel, Bounds on Degrees of Projective Schemes, Math.\ Ann.\ {\bf 302} (1995), 417--432.

\bibitem{Vi} R.\ H.\  Villarreal, {\it Monomial Algebras}, Dekker, New York, NY, 2001.



\end{thebibliography}
\end{document}